\documentclass[12pt]{amsart}
\usepackage{amsfonts, amsmath, amsthm}
\usepackage{psfig}

\newtheorem{pro}{Proposition}[section]
\newtheorem{thm}[pro]{Theorem}
\newtheorem{lem}[pro]{Lemma}

\newtheorem{cnj}[pro]{Conjecture}

\theoremstyle{definition}

\theoremstyle{remark}

\begin{document}

\title{Virtual homology of surgered torus bundles}
\author{Joseph D. Masters}

\begin{abstract}
Let M be a once-punctured torus bundle over $S^1$ with monodromy $h$.
  We show that, under certain hypotheses on $h$, ``most'' Dehn-fillings
 of $M$  (in some cases all but finitely many) are virtually
 $\mathbb{Z}$-representable.  We apply our results to show that
 even surgeries on the figure eight knot
 are virtually $\mathbb{Z}$-representable.
\end{abstract}
\maketitle

\section{Introduction}
Embedded incompressible surfaces are fundamental in the study
 of 3-manifolds. Accordingly, the following conjecture
 of Waldhausen and Thurston has attracted much attention:
\begin{cnj}
Let M be a closed, irreducible 3-manifold with infinite $\pi_1$.
  Then $M$ has a finite cover which is Haken.
\end{cnj}

The focus of this paper is the following, stronger, conjecture:
\begin{cnj}
Let $M$ be as above.  Then $M$ has a finite cover $\tilde{M}$ with
 $H_1(\tilde{M}, \mathbb{Z})$ infinite.
\end{cnj}

If $M$ is a compact 3-manifold, we say that $M$ is
 $\mathbb{Z}$-\textit{representable} if $H_1(M,\mathbb{Z})$ is infinite.
  If $M$ satisfies the conclusion of Conjecture 1.2, we say that $M$ is
 \textit{virtually} $\mathbb{Z}$\textit{-representable}.
  
  We shall give what appear to be the first examples of 3-manifolds with
 torus boundary for which all
 but finitely many fillings are virtually $\mathbb{Z}$-representable,
 but not $\mathbb{Z}$-representable (in fact non-Haken).  Boyer and Zhang
 have independently given examples of knot complements
 for which all but finitely many fillings are virtually Haken,
 but non-Haken [BZ].

  Before we can state our results, we must establish some notation.
  Let $F$ be a once-punctured torus with $\pi_1(F) = <[x], [y]>$, and
 basepoint $x_0 \in \partial F$
 (see Fig. 1).
\begin{figure}[htpb]
\begin{center}
\ \psfig{file=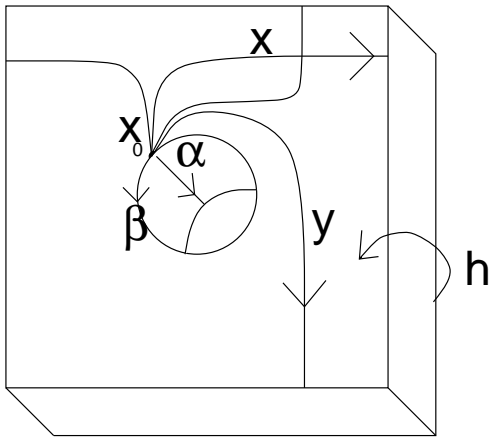}
\caption{Notation for the once-punctured torus bundle $M$.}
\end{center}
\label{notation}
\end{figure}
  Any orientation-preserving homeomorphism
 $h : F \rightarrow F$ is isotopic to one
 of the form $h = D_x^{r_1}D_y^{s_1}...D_x^{r_k}D_y^{s_k}$.
 Here $D_x$ and $D_y$ are Dehn twists along a simple closed curves homologous
 to $x$ and $y$, respectively.
 $D_x$ and $D_y$ induce the following actions on $\pi_1(F)$:\\
\begin{eqnarray*}
D_{x\sharp}(x) = x\\
D_{x\sharp}(y) = yx\\
D_{y\sharp}(x) = yx\\
D_{y\sharp}(y) = y
\end{eqnarray*}
  We may assume $h$ fixes $\partial F$.
  Let $M_h = (F \times I)/h$ be the once-punctured torus bundle
 with monodromy $h$.
  We specify a framing for $H_1(\partial M_h, \mathbb{Z})$ by setting the
 longitude $\beta = \partial F$ oriented counter-clockwise, and the meridian
 $\alpha = (x_0 \times I)/h$,
 where $x_0$ is some point on $\partial F$, and $\alpha$ is
 oriented as in Fig. 1.
  Then, for coprime integers $(\mu, \lambda)$, $M_h(\mu, \lambda)$
 denotes the manifold obtained by
 gluing a solid torus to
 $M_h$ in such a way that the curve $\alpha^{\mu}\beta^{\lambda}$
 becomes homotopically trivial.

We shall prove:
\begin{thm}\label{main}
Let $M_h$ be a once-punctured torus bundle over $S^1$, with
 monodromy $h = D_x^{r_1}D_y^{s_1}...D_x^{r_k}D_y^{s_k}$, and
 let $n = g.c.d \lbrace s_1, ..., s_k \rbrace$, $R = r_1 + ... + r_k$.
 Suppose $n$ is divisible by some $m$ such that\\
 i. $m \geq 6$ and $m$ is even or $m = 7$.  Then if $|\lambda| > 1$,
 all but finitely many
 Dehn-fillings $M_h(\mu, \lambda)$ are virtually $\mathbb{Z}$-representable.\\
 ii. $m \geq 5$, $m$ is odd, and $m \neq 7$.  Then if
 $1/|R\mu - \lambda| + 1/|R\mu - 2\lambda| +1/|\lambda| < 1$,
 $M_h(\mu, \lambda)$ is virtually $\mathbb{Z}$-representable. \\
 iii.  $m = 4$.  Then if $2/|R\mu - 2\lambda|+1/|\lambda| < 1$,
 $M_h(\mu, \lambda)$ is virtually $\mathbb{Z}$-representable.\\ 
\end{thm}

\textit{Remarks}: 1.  Analogous results hold if we replace
 $n$ by $gcd\lbrace r_1, ..., r_k \rbrace$ and $R$ by $s_1 + ... + s_k$.\\

  2. It was shown in [B1] that if $m \geq 2$, $n \geq 2$ and $mn \geq 8$
 but $mn \neq 9$,
 then all non-integral surgeries are virtually $\mathbb{Z}$-representable.
  In [B2] it was shown that if $4|n$, then for each $\mu$,
 $M_h(\mu, \lambda)$ is virtually $\mathbb{Z}$-representable for all
 but finitely many $\lambda$ coprime to $\mu$.\\

 3. From [CJR] and [FH], all but finitely many surgeries
 on a once-punctured torus bundle over $S^1$ yield non-Haken manifolds.\\

Theorem ~\ref{main} may be used to show that, for certain choices of $f$,
 all but finitely many surgeries on $M_f$ are
 virtually $\mathbb{Z}$-representable.
  For example:\\

\begin{thm} ~\label{fin}
Let $f = (D_xD_y)^{18}$.  Then every surgery on $M_f$ is virtually
 $\mathbb{Z}$-representable.
\end{thm}

The proof of Theorem ~\ref{fin} appears in Section 3.\\

In Section 3, we shall also prove the following theorems:\\
\begin{thm}\label{sister}
Let $(-1) = (D_xD_y^{-1}D_x)^2$, and let  
$N = M_{-D_xD_y}$ (also known as ``the figure eight knot's
 sister'').  Then if 
$1/|\mu - \lambda| + 1/|\mu - 2\lambda| + 1/|\lambda| < 1$,
 $N(\mu, \lambda)$ is virtually $\mathbb{Z}$-representable.
\end{thm}

\begin{thm}\label{fig8}
Let $K$ denote the figure-eight knot and let $M$ denote $S^3 - K$.
  Then, with respect to the canonical framing of knots in
 $S^3$,  any surgery of the form $M(2\mu, \lambda)$ is
 virtually $\mathbb{Z}$-representable.
 \end{thm}
  Other results on virtually $\mathbb{Z}$-representable figure-eight knot
 surgeries may be found in [M], [KL], [H], [N] and [B3].  In particular,
 it was shown in [KL] and [B3] that surgeries of the form
 $M(4\mu, \lambda)$ are virtually $\mathbb{Z}$-representable.
  It was also shown in [B3] that surgeries of the form
 $M(2\mu, \lambda)$ are virtually $\mathbb{Z}$-representable
 if $\lambda = \pm 7\mu$ (mod 15).

  Our techniques are extensions of Baker's.  The main new ingredient
 is the use of group theory to encode the combinatorics
 of cutting and pasting.

  I would like to thank Professor Alan Reid for his help and patience.

\section{Construction of covers}

  We begin by recalling Baker's construction of covering
 spaces of $M_h(\mu,\lambda)$ (see [B1], [B2]).
  Let $n$ be as in the statement of Theorem ~\ref{main}, and let
 $\hat{F}$ be the  $kn$-fold cover of $F$ associated to the kernel
 of the map $\phi : \pi_1(F) \rightarrow \mathbb{Z}_k \times \mathbb{Z}_n$,
 with $\phi([x]) = (1,0)$ and $\phi([y]) = (0,1)$ (see Fig. 2).
  
\begin{figure}[htbp]

\begin{center}

\ \psfig{file=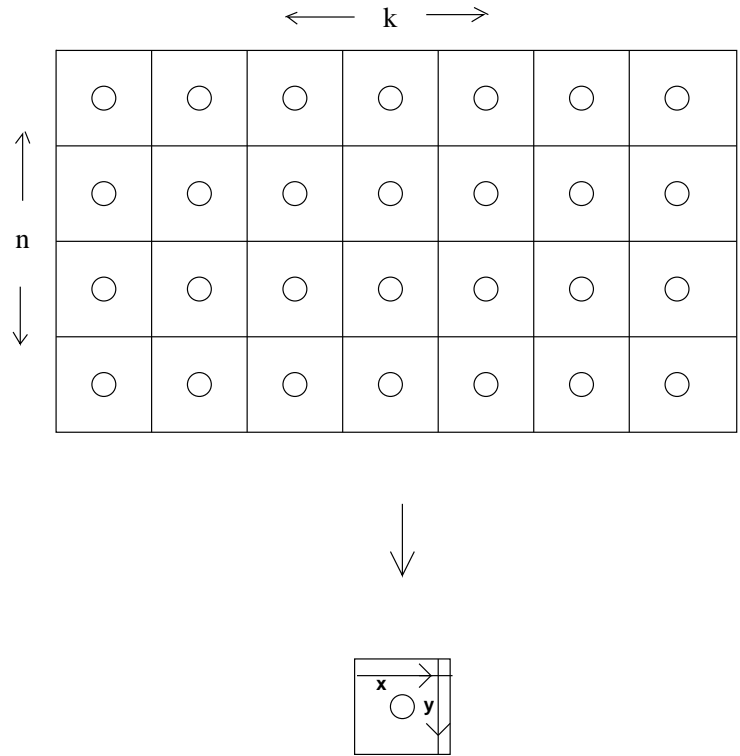}

\caption{The cover $\hat{F}$ of $F$}

\end{center}
\label{coverbasic}
\end{figure}

Now create a new cover, $\tilde{F}$, of $F$ by making vertical cuts in
 each row of $\hat{F}$,
 and gluing the left side of each cut to the right side of another cut in
 the same row.  An example is pictured in Figure 3,
 where the numbers
 in each row indicate how the edges are glued.

\begin{figure}[htbp]

\begin{center}

\ \psfig{file=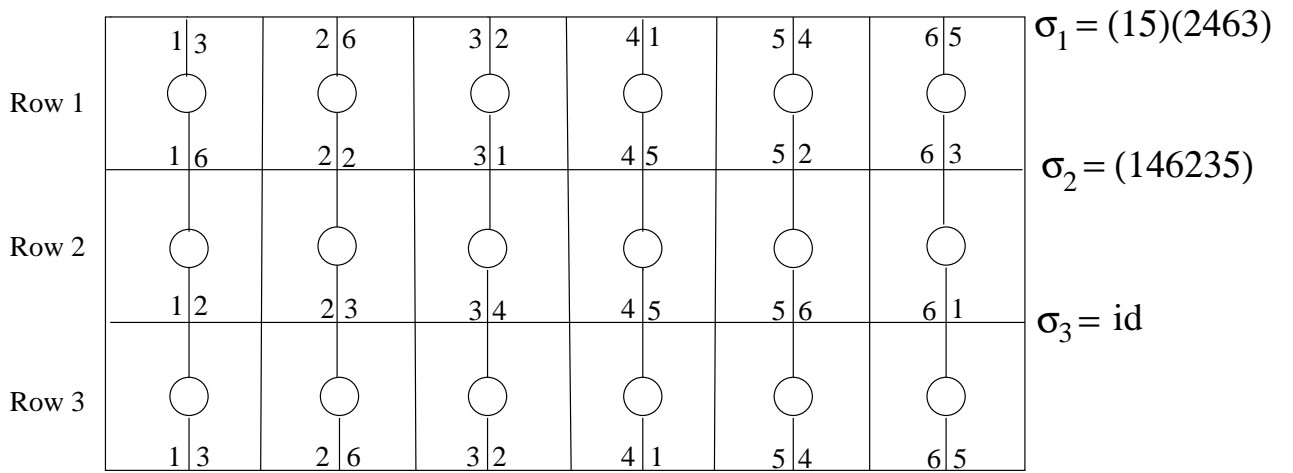}

\caption{The permutations encode the combinatorics of the gluing}

\end{center}
\label{comb}
\end{figure}

  If $h$ lifts to a map $\tilde{h} : \tilde{F} \rightarrow \tilde{F}$,
 then the mapping cylinder $\tilde{M_h} = \tilde{F}/\tilde{h}$ is a
 cover of $M_h$.
  Furthermore, if the loop $\alpha^{\mu}\beta^{\lambda}$ lifts to loops in
 $\tilde{M_h}$, then
 the cover extends to a cover $\tilde{M_h}(\mu,\lambda)$ of
 $M_h(\mu, \lambda)$.

If the cover $\tilde{M}_h$ exists, then we may estimate its homology with
 the formula $b_1(\tilde{M}_h) = rank(fix(\tilde{h}_*))$, where $\tilde{h}_*$
 is the map on $H_1(\tilde{M}, \mathbb{Z})$ induced by $\tilde{h}$, and
 $fix(\tilde{h}_*)$ is the subgroup of $H_1(\tilde{M}, \mathbb{Z})$ fixed
 by $\tilde{h}_*$ (see [H] for a proof).  We shall use this
 formula to prove that, in some cases,
 $b_1(\tilde{M}) >$ \textit{number of boundary components of} $\tilde{M}$,
 which ensures that $b_1(\tilde{M}(\mu, \lambda)) > 0$.

We now introduce some notation to describe the cuts of $\tilde{F}$
 (see Fig. 3). 
 $\tilde{F}$ is naturally divided into rows, which we label
 1, ..., n.  The cuts divide each row
 into pieces, which we label $1, ..., k$.  If we slide a point in the top
 half of
 the $i^{th}$ row through the cut to its right, we induce a permutation on
 $\lbrace 1, ...,k \rbrace$, which we denote $\sigma_i$.  Thus the cuts on
 $\tilde{F}$
 may be encoded by elements $\sigma_1, ..., \sigma_n \in S_k$, the
 permutation group on $k$ letters.

Next, we find algebraic conditions on the $\sigma_i$'s which will guarantee
 that the cover of $F$ extends to a cover of $M(\mu, \lambda)$.  We first
 must pick $k, n,$ and $\lbrace \sigma_1, ..., \sigma_n \rbrace$ so that
 $h$ lifts to $\tilde{F}$.

\begin{lem} \label{Fcover}
If \\
I. $[\sigma_i, \sigma_1\sigma_2...\sigma_{i-1}] = 1\,$ for all  $i$ and\\
II. $\sigma_1\sigma_2...\sigma_n = 1$\\
then $h$ lifts to $\tilde{F}$.
\end{lem}

\begin{proof}

Note that $D_y^n$ lifts to Dehn twists on $\tilde{F}$.
  Therefore, 
 we need only ensure that $D_x$ lifts.  We shall attempt to lift $D_x$
 to a sequence of ``fractional Dehn twists'' along the rows of $\tilde{F}$.
  Let $\tilde{x}_i$ denote the disjoint union of the lifts of $x$
 to the $i^{th}$ row of $\tilde{F}$.
  We first attempt to lift $D_x$ to row 1, twisting $1/k^{th}$ of the way
 along $\tilde{x}_1$.
  Considering the action on the bottom half of row 1, we find that
 the cuts
 are now matched up according to the permutation
 $\sigma_1^{-1}\sigma_2\sigma_1$.  Thus, for $D_x$ to lift to row 1
 we assume $\sigma_1$ and $\sigma_2$ commute.  We now twist along
 $\tilde{x}_2$.  The top halves of the squares in row 2 are moved
 according to the permutation
 $\sigma_1\sigma_2$, and the lift will extend to all of row 2 if and only if
 $\sigma_3$ commutes with $\sigma_1\sigma_2$.  We continue in this manner,
 obtaining the conditions in I.  After we twist through
 $\tilde{x}_n$, we need to be back where we started in row 1, so we require
 the additional condition $\sigma_1\sigma_2...\sigma_n = 1$.

\end{proof}

Note that the loop $\alpha^{\mu}$ lifts to loops in $\tilde{M_h}$
 if $\tilde{h}^{\mu} = id$, and that
 the loop $\beta^{\lambda}$ lifts to loops in $\tilde{M_h}$
 if $\sigma_{i+1}\sigma_i^{-1} = id$ for all $i = 1, ..., n$.
 Then, by considering the action of $\tilde{h}$
 on $\tilde{M_h}$,
 the following condition for a loop in $\partial M_h$ to lift to
 $\tilde{M_h}$ is easily verified:

\begin{lem} \label{Mcover}
The loop $\alpha^{\mu}\beta^{\lambda} \subset \partial M_h$ lifts to
 loops in
 $\tilde{M_h}$ if and only if\\ 
III. $(\sigma_1...\sigma_i)^{R\mu}(\sigma_{i+1}\sigma_i^{-1})^{\lambda}=1$,
 for i= 1, ..., n.
\end{lem}

Therefore we may construct covers of $M_h(\mu, \lambda)$ simply by finding
 permutations satisfying conditions I-III.

\vskip1pc

\begin{proof} (Of Theorem ~\ref{main})\\

\textbf{Case 1:} m = 4.

\vskip1pc

\textit{Construction of the cover of $M_h(\mu, \lambda)$:}

To construct a cover of $M_h(\mu, \lambda)$, we must first construct a cover
 of $F$.  It was shown in the discussion prior to Lemma \ref{Fcover}
 that there is a unique such cover associated to any four permutations
 $\sigma_1,\sigma_2, \sigma_3$ and $\sigma_4$ in any permutation
 group $S_k$.

 To ensure that the cover of $F$ extends to a cover of $M_h$,
 we shall set $\sigma_2 = \sigma_1^{-1}$
 and $\sigma_4 = \sigma_3^{-1}$ (see Fig. 4a).  Then conditions
 I and II of Lemma \ref{Fcover} are satisfied automatically,
 so that any choice of
 $\sigma_1$ and $\sigma_3$ will determine a cover of $M_h$.

To ensure that the cover extends to $M_h(\mu, \lambda)$, we
 must arrange for the surgery curve $\alpha^{\mu}\beta^{\lambda}$ 
 to lift to $\tilde{M_h}$.  By Lemma \ref{Mcover}, $\alpha^{\mu}\beta^{\lambda}$
 will lift provided that 
 $\sigma_1, ..., \sigma_4$ satisfy condition III, which reduces to:

\begin{eqnarray} \label{cover4eqns}
\sigma_1^{R\mu-2\lambda} = 1\\
(\sigma_3\sigma_1)^{\lambda} = 1\\
\sigma_3^{R\mu-2\lambda}=1\\
(\sigma_1\sigma_3)^{\lambda} = 1
\end{eqnarray}

\begin{figure}[htpb] \label{cover45}

\begin{center}

\ \psfig{file=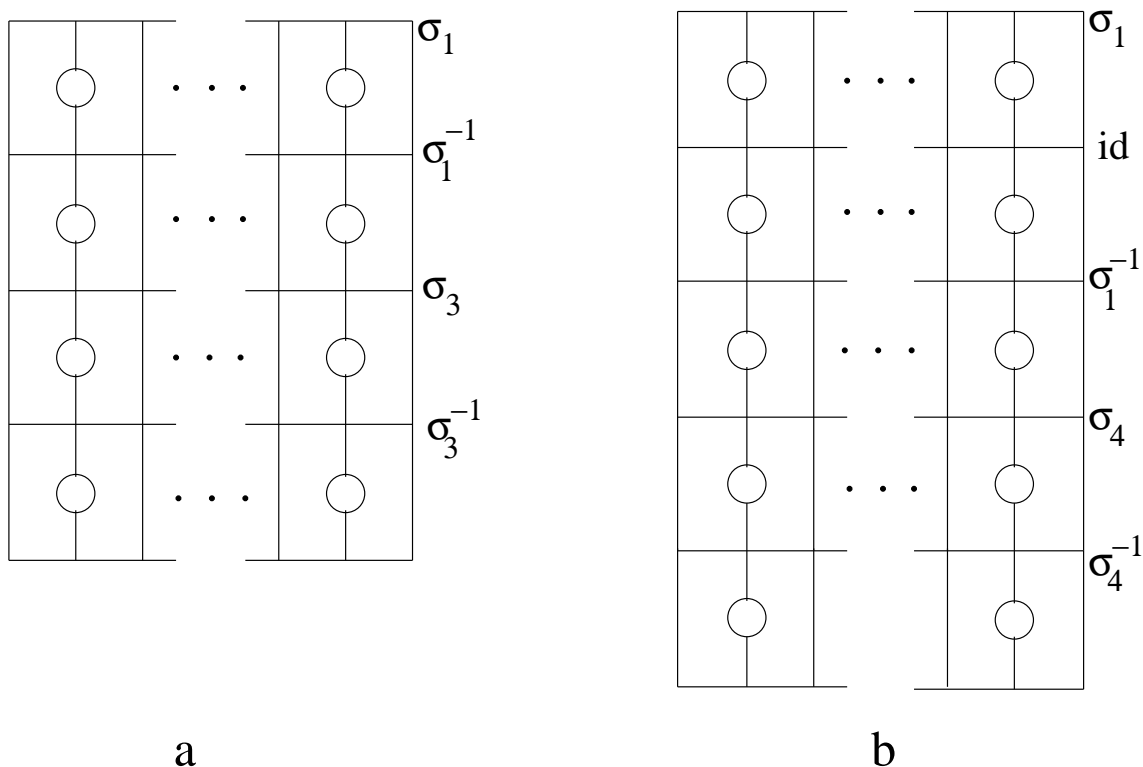}

\caption{a. The cover when $n = 4.\,\,\,\,\,\,\,\,$ b. The cover when $n = 5$.}

\end{center}

\end{figure}

 Any pair of permutations $\sigma_1$ and $\sigma_3$ satisfying
 equations (1)-(4) determines a unique cover of $M_h(\mu, \lambda)$.
  We now turn our attention to the construction of such permutations.

 Consider the abstract group $G$
 generated by the symbols $\bar{\sigma}_1$ and $\bar{\sigma}_3$,
 satisfying relations (1)-(4).  $G$ is a
 $(|R\mu-2\lambda|, |R\mu-2\lambda|, |\lambda|)$-triangle group.
 It is well-known that if
 $1/|R\mu - 2\lambda| + 1/|R\mu-2\lambda|+ 1/|\lambda| < 1$,
 then $G$ is residually finite, and hence surjects a finite group $H$ such that
 the images of $\bar{\sigma}_1$, $\bar{\sigma}_3$,
 and $\bar{\sigma_3}\bar{\sigma_1}$ have order
 $|R\mu - 2\lambda|$.
 By taking the permutation representation of $H$, we then obtain
 permutations $\sigma_1$ and $\sigma_3$ satisfying conditions 
 (1)-(4).  Note that the permutations act on $|H|$ letters, so $\tilde{M}$
 is a $4|H|$-fold cover of $M_h$.

 Associated with the permutations $\sigma_1$
 and $\sigma_3$ we have covers $\tilde{M_h}$ and $\tilde{F}$ of $M_h$ and $F$,
  and a cover $\tilde{M_h}(\mu, \lambda)$ of $M_h(\mu, \lambda)$;
 \\
\\
\textit{ Claim}: $b_1(\tilde{M_h}(\mu, \lambda)) > 0$.\\
\\
\textit{Proof of claim:} It suffices to show that $\tilde{h}_*$ has 
  a non-peripheral class $[\delta] \in H_1(\tilde{F})$ with
 $\tilde{h}_{\ast}([\delta]) = [\delta]$.
 To construct this element, we shall first find a non-peripheral class
 $[\delta_2]$ in row 2, as follows.

\begin{figure}[htbp] \label{F2}
\begin{center}
\ \psfig{file=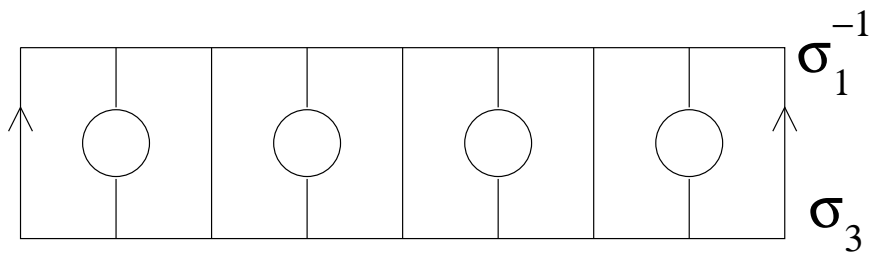}
\caption{The surface $\tilde{F}_2$ (with $|H| = 4$).}
\end{center}
\end{figure}

  Consider the sub-surface $\tilde{F}_2$ obtained
 by deleting rows 1, 3 and 4 from
 $\tilde{F}$ (see Fig. 5).  The punctures of $\tilde{F}_2$
 are in 1-1 correspondence with the cycles of $\sigma_1$, $\sigma_3$ and
 $\sigma_3\sigma_1$.  Any permutation $\tau$ coming from the permutation
 representation of $H$ decomposes as a product of
 $|H|/order(\tau)$ disjoint $order(\tau)$-cycles.
  Therefore $\tilde{F}_2$ has
 $|H|(1/order(\sigma_1)+ 1/order(\sigma_3)+ 1/order(\sigma_3\sigma_1)) < |H|$
 punctures.  Since $\chi(\tilde{F}_2) = -|H|$, we deduce that
 $\tilde{F}_2$ contains a non-peripheral class $[\delta_2]$.
 $\delta_2$ also represents a non-peripheral class
 in $\tilde{F}$, since it has non-zero intersection
 number with a class of $\tilde{F}$ in row 2.

\begin{figure}[htbp] \label{cancel}
\begin{center}
\ \psfig{file=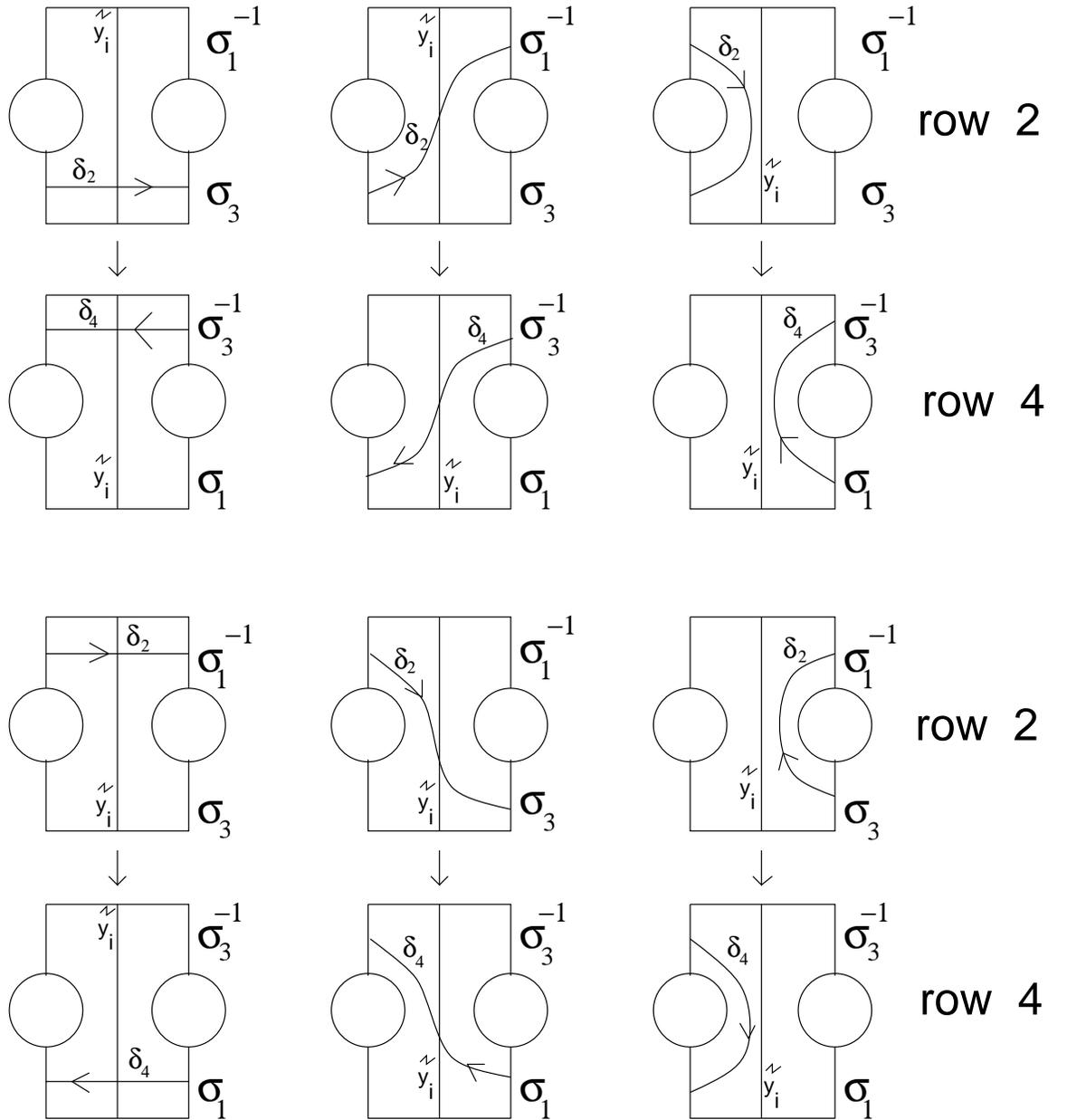}
\caption{How to find canceling loops in rows 2 and 4}
\end{center}
\end{figure}

  We may find a corresponding non-peripheral loop
 $\delta_4$ in row 4, such that
 $I([\delta_2+\delta_4, [\tilde{y}_i]) = 0 \,$ for all $i$ (see
 Fig. 6 for
 the notation and the idea of the proof).
 Let $[\delta] = [\delta_2 + \delta_4]$.  Then, since $[\delta]$ has
 non-zero intersection number with classes in row 2 and row 4,
 it is a non-peripheral class.  Since $I[\delta, \tilde{y}_i] =0\,$
 for all $i$,
 it is fixed by $\tilde{D_y^4}_*$, and since $\tilde{D}_x$ fixes rows
 2 and 4, it is fixed by $\tilde{D_x}_*$.  Therefore it is
 fixed by $\tilde{h}_*$,
 concluding the proof of the claim, and of Case 1.

\vskip1pc

\textbf{Case 2:} $m \geq 5$ and $m$ is odd.

\vskip1pc

\textit{Case 2a:} $m = 5$.\\
The construction proceeds analogously to the case m = 4.
  We require permutations $\sigma_1, ..., \sigma_5$
 satisfying conditions I - III.  Again, to simplify matters,
 we shall impose some extra conditions:
 $\sigma_2 =\, id$, $\sigma_3 = \sigma_1^{-1}$, and $\sigma_5 = \sigma_4^{-1}$
  (see Fig. 4).  Then conditions I-III reduce to:
\begin{eqnarray*}
\sigma_1^{R\mu-\lambda} = 1\\
(\sigma_1\sigma_4)^{\lambda} = (\sigma_4\sigma_1)^{\lambda} = 1\\
\sigma_4^{R\mu-2\lambda} = 1
\end{eqnarray*}
Again, these relations determine a triangle group, which, under the hypotheses
 on $\mu$ and $\lambda$, is hyperbolic.  The rest of the proof is identical
 to Case 1, except that now the fixed class is in rows 3 and 5.

\vskip1pc

\textit{Case 2b:} $m \geq 9$ and $m$ is odd.

\vskip1pc

Consider the cover obtained by setting $\sigma_2 = \sigma_1^{-1}$,
 $\sigma_4 = id$, $\sigma_5 = \sigma_3^{-1}$,
 $\sigma_6 = \sigma_1$, $\sigma_7 = \sigma_1^{-1}$,
  and for $i = 4, ..., k$, $\sigma_{2i+1} = \sigma_{2i}^{-1}$
 (see Fig. 7a).

\begin{figure}[htbp] \label{cover89}

\begin{center}

\ \psfig{file=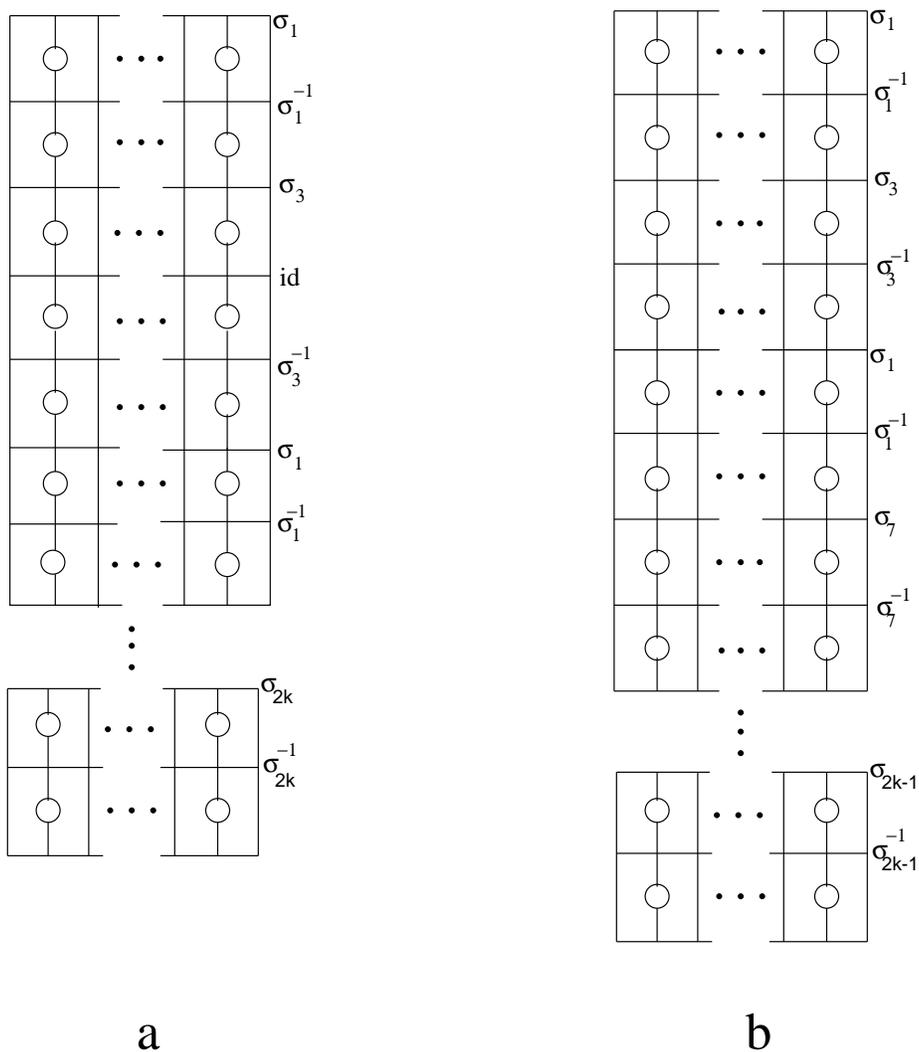}

\caption{a. The cover for $n= 2k + 1 \geq 9$. $\,\,\,$ b. The cover for
 $n = 2k \geq 8$.}

\end{center}

\end{figure}
 
 The corresponding relations are:
\begin{eqnarray} \label{Coxeqns}
\sigma_1^{R\mu-2\lambda} = 1 \,\,\,\,\,\,\,\,\,\,\,\,\,\,\,\,\,\,\,\,\,\,\,\,\,\,\,\,\,\,\,\,\,\,\,\,\,\,\,\,\,\,\,\,\,\,\,\,\\
\sigma_3^{R\mu -\lambda} = 1 \,\,\,\,\,\,\,\,\,\,\,\,\,\,\,\,\,\,\,\,\,\,\,\,\,\,\,\,\,\,\,\,\,\,\,\,\,\,\,\,\,\,\,\,\,\,\,\,\\
(\sigma_3\sigma_1)^{\lambda} = (\sigma_1\sigma_3)^{\lambda} = 1 \,\,\,\,\,\,\,\,\,\,\,\,\,\,\,\,\,\,\,\,\,\,\,\,\,\,\,\,\,\,\\
(\sigma_8\sigma_1)^{\lambda} = 1 \,\,\,\,\,\,\,\,\,\,\,\,\,\,\,\,\,\,\,\,\,\,\,\,\,\,\,\,\,\,\,\,\,\,\,\,\,\,\,\,\,\,\,\,\,\,\,\,\\
\sigma_{2i}^{R\mu - 2\lambda} = 1 \textnormal{ for } i = 4, ..., k \,\,\,\,\,\,\,\,\,\, \\
(\sigma_{2i+2}\sigma_{2i})^{\lambda} = 1 \textnormal{ for } i = 4, ..., k-1\\
(\sigma_{1}\sigma_{2k})^{\lambda} = 1 \,\,\,\,\,\,\,\,\,\,\,\,\,\,\,\,\,\,\,\,\,\,\,\,\,\,\,\,\,\,\,\,\,\,\,\,\,\,\,\,\,\,\,\,\,\,\,\,
\end{eqnarray}

These relations again determine a Coxeter group.  It is well-known
 (see [V]) that any such group surjects a finite group ``without
 collapsing''-- i.e such that the orders of the images of the
 $\sigma_i$'s and $\sigma_i\sigma_j$'s are as given in (5)-(11).
  Then, arguing as in Case 1,  we may find a non-peripheral fixed class in rows
 2 and 5.
\vskip1pc
\textbf{Case 3:} n = 6
\vskip1pc

\textit{Case 3a:} $2/|R\mu - \lambda| + 1/|\lambda| < 1$

Again, we need permutations $\sigma_1, ..., \sigma_6$ satisfying I - III.
  In this case we impose the additional conditions $\sigma_2 = id$,
 $\sigma_3 = \sigma_1^{-1}$, $\sigma_5 = id$, and $\sigma_6 = \sigma_4^{-1}$.
  Then conditions I-III reduce to:
\begin{eqnarray*}
\sigma_1^{R\mu-\lambda} = 1\\
(\sigma_1\sigma_4)^{\lambda} = (\sigma_4\sigma_1)^{\lambda} = 1\\
\sigma_4^{R\mu - \lambda} = 1
\end{eqnarray*}

These relations determine a triangle group,
 and we find a fixed class in rows 3 and 6.\\
\\

\textit{Case 3b:}  $|\lambda| > 2$ and $|R\mu - 3\lambda| \geq |\lambda|$, or
 $\lambda$ is even (non-zero) and $|R\mu - 3\lambda| \geq 4$.

When $n = 3$, conditions I-III may be abelianized to obtain a cyclic
 group of order $|R\mu + 3\lambda|$.  Specifically,
 they are satisfied by setting
 $\sigma_1 = (1, \, 2, \,...\,  ,R\mu - 3\lambda)$,
 $\sigma_2 = \sigma_1^{-2}$, and $\sigma_3 = \sigma_1$.
 For $n = 6$, we may ``double'' this cover:
 that is take $\sigma_1, \sigma_2, \sigma_3$ as above, and then
 set $\sigma_4 = \sigma_1$, $\sigma_5 = \sigma_2$, and $\sigma_6 = \sigma_3$.
  Then we modify the corresponding cover
 $\tilde{M}(\mu,\lambda)$ of $M(\mu, \lambda)$
 by making  horizontal cuts in adjacent squares of row 3 and
 gluing the flaps back together as indicated by Fig. 8.
 If $\lambda$ is even, we make two non-adjacent cuts and glue the top of
 one to the bottom of the other.
 If $\lambda$ is odd, we make $(|\lambda| - 1)$/2 pairs of adjacent cuts and
 glue the top of the one cut to the bottom of the other cut in its pair.
  Now make the same  
 cuts in row 6, with the same identifications.  Since rows 3 and 6
 are fixed by $\tilde{D}_x$, $D_x$ still lifts to the modified
 $\tilde{M_h}(\mu, \lambda)$, and since the $\tilde{y}$'s still project 6 to 1
 onto $y$, $D_y$ lifts also;  so $h$ lifts.  Also, one may check
 that $\alpha^{\mu}\beta^{\lambda}$ still lifts, so $\tilde{M_h}(\mu, \lambda)$
 remains a cover of $M_h(\mu, \lambda)$.

  To see that $b_1(\tilde{M_h}(\mu, \lambda)) > 0$, note that $\tilde{D}_x$
 fixes rows 3 and 6, so it is enough to find a non-peripheral loop in
 row 3 and add it to the corresponding loop in row 6 with opposite
 orientation.  As in Case 1, the existence of such a non-peripheral loop
 follows from an Euler characteristic argument (or see Fig. 8).

\vskip1pc
Note that Case 3a or 3b applies to all but finitely many $(\mu, \lambda)$
 with $|\lambda| > 1$.

\begin{figure}[htbp] \label{cover6}

\begin{center}

\ \psfig{file=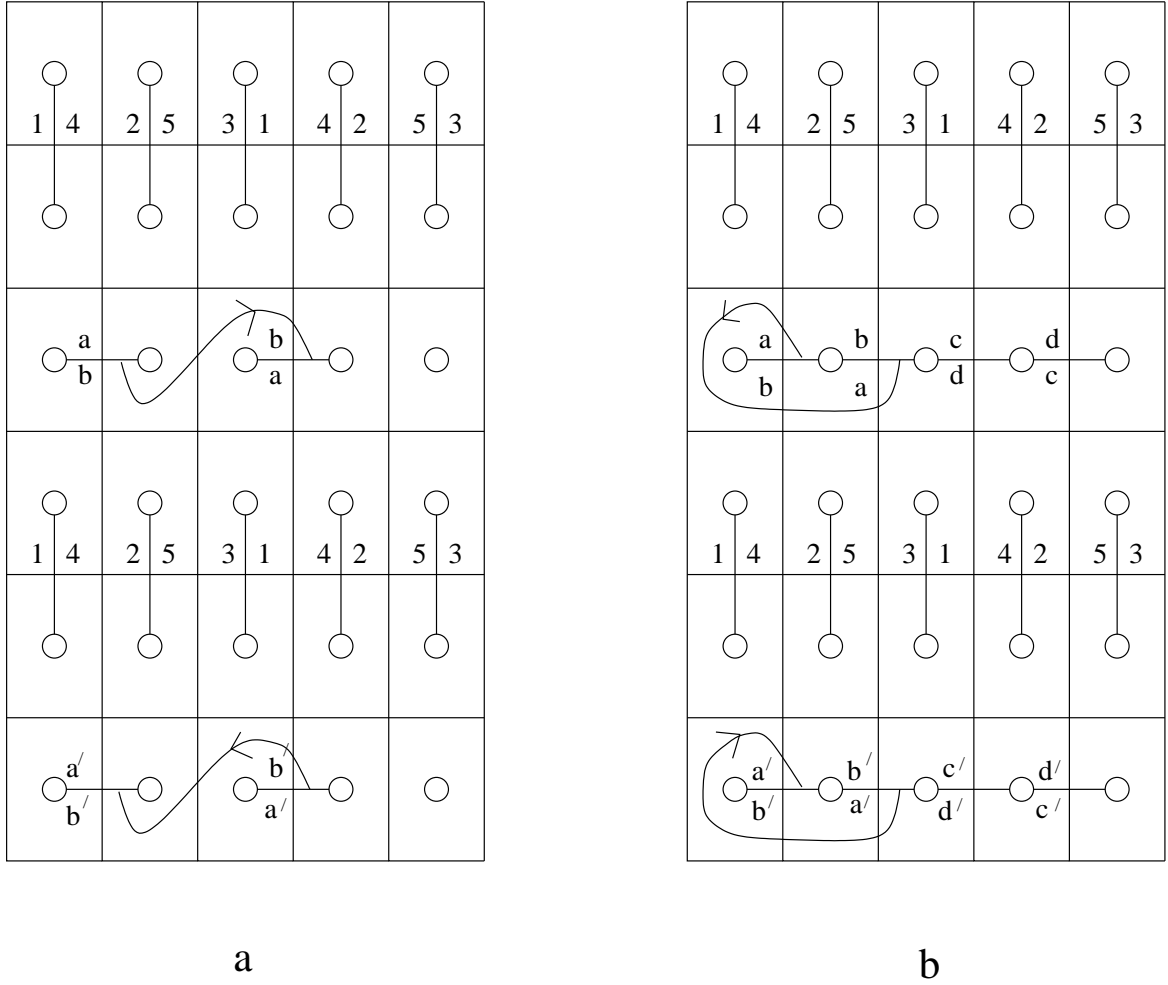}

\caption{a. The cover and fixed class for $n=6, R\mu = 1, \,\, \lambda=2$
 b. The cover and fixed class for $n = 6, R\mu = 10, \,\, \lambda = 5$}

\end{center}

\end{figure}

\vskip1pc
\textbf{Case 4:} $n = 2k \geq 8$.
\vskip1pc

\textit{Case 4a:} $2/|R\mu-2\lambda| + 1/|\lambda| < 1$.
Set $\sigma_2 = \sigma_1^{-1}$, $\sigma_4 = \sigma_3^{-1}$,
 $\sigma_5 = \sigma_1$, $\sigma_6 = \sigma_1^{-1}$, and
 $\sigma_{2i} = \sigma_{2i-1}^{-1}$ for $i = 4, ..., k$ (see Fig. 
 7b).
  Then, as in Case 2, these relations determine a Coxeter group.
  We may find a non-peripheral fixed class in rows 2 and 4.

\vskip1pc

\textit{Case 4b:} $|R\mu - \lambda| \leq 2$
We cannot guarantee, in this case, that there will always be
 a cover with $b_1 > 0$, but we shall show that there are
 at most finitely many exceptions. 

We argue as in Case 3b.  
Take permutations $\sigma_1, ..., \sigma_k$, and consider the relations
 obtained by abelianizing conditions I-III.  We claim that they can be
 satisfied by setting $\sigma_1 = (1, \,2, \,3, \, ..., \, N)$, for some N, and
 setting each $\sigma_i$ to an appropriate power of $\sigma_1$.
  We have already seen that this may be done when $k = 3$.

The $\sigma_i$'s must satisfy the following conditions:
\begin{equation} \label{sighardeqn1}
\sigma_1^{R\mu-\lambda}\sigma_2^{\lambda} = 1
\end{equation}
\begin{equation} \label{sighardeqn2}
\sigma_1^{R\mu}\sigma_2^{R\mu-\lambda}\sigma_3^{\lambda} = 1
\end{equation}

\begin{center}
$.\atop{.\atop.}$
\end{center}
\begin{equation} \label{sighardeqn3}
\sigma_1^{R\mu}\sigma_2^{R\mu}...\sigma_{k-2}^{R\mu}
\end{equation}
\begin{equation} \label{sighardeqn4}
\sigma_{k-1}^{R\mu-\lambda}\sigma_k^{\lambda}=1
\end{equation}
\begin{equation} \label{sighardeqn5}
\sigma_1^{R\mu+\lambda}\sigma_2^{R\mu}...\sigma_{k-1}^{R\mu}
\end{equation}
\begin{equation} \label{sighardeqn6}
\sigma_k^{R\mu-\lambda} = 1
\end{equation}
\begin{equation} \label{sighardeqn7}
\sigma_1\sigma_2...\sigma_k = 1
\end{equation}

We shall assume that this system has a cyclic solution,
 so we may substitute $\sigma_i = \sigma_1^{e_i}$.
Then, equations (~\ref{sighardeqn1})-(~\ref{sighardeqn7})
 are equivalent to the following conditions
 on the exponents (all of the following equations in this case are
 taken mod N):
\begin{equation} \label{exphardeqn1}
R\mu-\lambda+\lambda e_2 = 0
\end{equation}
\begin{equation} \label{exphardeqn2}
R\mu + (R\mu-\lambda)e_2 +\lambda e_3 = 0
\end{equation}
$.\atop{.\atop{.}}$
\begin{equation}\label{exphardeqn3}
R\mu+R\mu e_2+...+R\mu e_{k-2}+(R\mu-\lambda)e_{k-1}+\lambda e_k= 0
\end{equation}
\begin{equation}\label{exphardeqn4}
R\mu+\lambda+R\mu e_2+...+R\mu e_{k-1}+(R\mu-\lambda)e_k=0
\end{equation}
\begin{equation}\label{exphardeqn5}
1+e_2+...+e_k=0
\end{equation}
  (~\ref{exphardeqn4}) and (~\ref{exphardeqn5})
 imply that $\lambda = \lambda e_k$.  Let us
 set $e_k = 1$, eliminating equation (~\ref{exphardeqn4}).
  Then, using (~\ref{exphardeqn5}), we may pair off
 (~\ref{exphardeqn1}) and (~\ref{exphardeqn3}) to deduce that
 $\lambda e_2 = \lambda e_{k-1}$, and we set $e_2 = e_{k-1}$ to
 eliminate (~\ref{exphardeqn3}).
  Similarly, we set $e_3 = e_{k-2}$, and so on.  If $k$ is even, we are left
 with equations:

\begin{equation} \label{redhardeqn1}
R\mu-\lambda + \lambda e_2 = 0
\end{equation}
\begin{equation}\label{redhardeqn2}
R\mu + (R\mu - \lambda) e_2 + \lambda e_3 = 0
\end{equation}
\begin{center}
$.\atop{.\atop.}$
\end{center}
\begin{equation} \label{redhardeqn3}
R\mu+R\mu e_2 + ... + (R\mu-\lambda)e_{k/2-1} + \lambda e_{k/2} = 0
\end{equation}
\begin{equation} \label{redhardeqn4}
R\mu+R\mu e_2 + ... + (R\mu-\lambda)e_{k/2} + \lambda e_{k/2} = 0
\end{equation}
\begin{equation} \label{redhardeqn5}
2 + 2e_2 + ... + 2e_{k/2} = 0
\end{equation}
If we replace (~\ref{redhardeqn5}) by
\begin{equation}\label{redhardeqn5'}
1 + e_2 + ... + e_{k/2} = 0
\end{equation} 
  then we may eliminate (~\ref{redhardeqn4}).  Then solve
 for $\lambda e_2, \lambda^2 e_3, ..., \lambda^{k/2-1}e_{k/2}$.
  By (~\ref{redhardeqn5'}), we have:
\begin{eqnarray}
\lambda^{k/2-1} + \lambda^{k/2-2}(\lambda e_2)
 + \lambda^{k/2-3}(\lambda^2e_3)+
 ... + \lambda^{k/2-1}e_{k/2}=0
\end{eqnarray}
  Substituting our solutions for $\lambda e_2, \lambda^2 e_3$ and so on,
 we get the equation $N = 0$ for some integer $N$; the system
 has a solution in $\mathbb{Z}/N\mathbb{Z}$.

 If $k$ is odd, then our reduced system looks like:
\begin{equation}\label{redoddeqn1}
R\mu-\lambda + \lambda e_2 = 0
\end{equation}
\begin{equation}\label{redoddeqn2}
R\mu + (R\mu - \lambda) e_2 + \lambda e_3 = 0
\end{equation}

\begin{center}
$.\atop{.\atop.}$
\end{center}
\begin{equation}\label{redoddeqn3}
R\mu+R\mu e_2 + ... + (R\mu-\lambda)e_{(k-1)/2} + \lambda e_{(k+1)/2} = 0
\end{equation}
\begin{equation}\label{redoddeqn4}
R\mu+R\mu e_2 + ... + (R\mu-\lambda)e_{(k+1)/2} + \lambda e_{(k-1)/2} = 0
\end{equation}
\begin{equation}\label{redoddeqn5}
2 + 2e_2 + ... + 2e_{(k-1)/2} + e_{(k+1)/2} = 0
\end{equation}
 Adding (~\ref{redoddeqn3}) and (~\ref{redoddeqn4}) gives a multiple of
 (~\ref{redoddeqn5}), so we may eliminate (~\ref{redoddeqn4}).
  Then we solve
 for $\lambda e_2, \lambda^2 e_3, ..., \lambda^{(k-1)/2}e_{(k+1)/2}$.
  By (~\ref{redoddeqn5}), we have:
\begin{eqnarray*}
2\lambda^{(k-1)/2} + 2\lambda^{(k-3)/2}(\lambda e_2)
 + 2\lambda^{k-5/2}(\lambda^2e_3)+ ...\\
 + 2\lambda(\lambda^{(k-3)/2}e_{(k-1)/2}) + \lambda^{(k-1)/2}e_{(k+1)/2}=0
\end{eqnarray*}
 And again we get a solution in $\mathbb{Z}/N\mathbb{Z}$ for some $N$.

Then, as in Case 3b, $M(\mu, \lambda)$ will have a cover with $b_1 > 0$,
 provided that $|N| \geq |\lambda|$ and $|\lambda| > 2$.
  Solving for $N$, if $k$ is even, gives:
\begin{equation}\label{Neveneqn}
 N = \lambda^{k/2-1} + \lambda^{k/2-2}(\lambda-R\mu)
 + \lambda^{k/2-3}[(\lambda- R\mu)^2 - R\mu\lambda]
\end{equation}
\begin{eqnarray*}
   + \lambda^{k/2-4}[(\lambda -R\mu)((\lambda -R\mu)^2-R\mu\lambda)-R\mu\lambda(\lambda -R\mu) -R\mu\lambda^2]
 + ...
\end{eqnarray*}
  and if $k$ is odd:
\begin{equation}\label{Noddeqn}
N = 2\lambda^{(k-1)/2} + 2\lambda^{(k-3)/2}(\lambda -R\mu)
+ 2\lambda^{(k-5)/2}[(\lambda -R\mu)^2 - R\mu\lambda]
\end{equation}
\begin{eqnarray*}
  & & + 2\lambda^{(k-7)/2}[(\lambda - R\mu)((\lambda -R\mu)^2-R\mu\lambda)-R\mu\lambda(\lambda -R\mu) -R\mu\lambda^2]\\
  & & + ... + 1[..] 
\end{eqnarray*}
  We are supposing that $|R\mu -\lambda| \leq 2$.  So for large
 $\mu$ or $\lambda$, $R\mu/\lambda \rightarrow 1$,
 and for $k$ even,
\begin{eqnarray*}
N = o[\lambda^{k/2-1} + \lambda^{k/2-3}(-\lambda^2)
 + \lambda^{k/2-4}(-\lambda^3) + ...]\\
 = o[\lambda^{k/2-1} - \lambda^{k/2-1} - \lambda^{k/2-1} - ...]
\end{eqnarray*}
 So if $k$ is even and $k\geq 8$, then for all but
 finitely many $\mu$ and $\lambda$, $|N| > |\lambda|$, and we are done.
 Similarly, if $k$ is odd and $k \geq 7$, then we are done.
 In the remaining cases, $|N|$ is given by:
\begin{eqnarray*}
 k = 4,\,\,\, |N| = |R\mu - 2\lambda|\\
 k = 5,\,\,\, |N| = |(R\mu)^2 - 5R\mu\lambda + 5\lambda^2|\\
 k = 6,\,\,\, |N| = |(R\mu)^2 - 4R\mu\lambda + 3\lambda^2|
\end{eqnarray*}
  One may check that each condition is satisfied by only finitely many
 relatively prime pairs $(\mu, \lambda)$ with $|R\mu -\lambda| \leq 2$.
 This concludes the proof in Case 4b.

\vskip1pc

\textbf{Case 5:} $n = 7$, and $|\lambda| > 1$.
\vskip1pc

\textit{Case 5a:} $1/|R\mu -\lambda| + 1/|\lambda| < 2/3$
 and $|(R\mu-2\lambda)^2 - 2\lambda^2| > 2|R|$.

\vskip1pc

We shall consider covers with $\sigma_2 = id$, $\sigma_3=\sigma_1^{-1}$,
$\sigma_6 = \sigma_5^{-1}$, and $\sigma_7 = \sigma_4^{-1}$ (see Fig.
 9a).

\begin{figure}[htbp] \label{7}
\begin{center}
\ \psfig{file=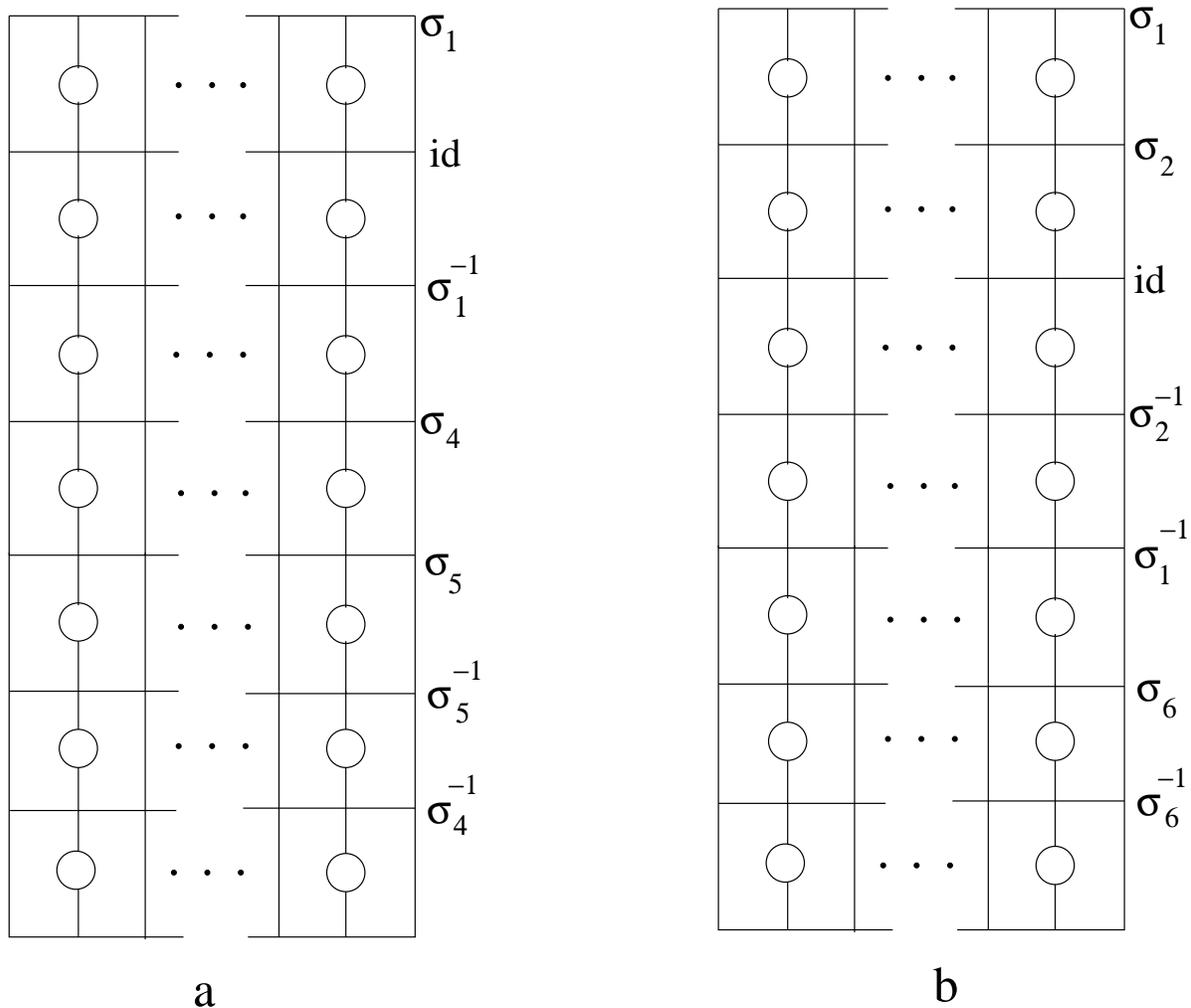}
\caption{Two covers for $n = 7$.}
\end{center}
\end{figure}

  We obtain conditions:
\begin{equation}\label{cover7eqn1}
[\sigma_4,\sigma_5] = 1
\end{equation}
\begin{equation}\label{cover7eqn2}
\sigma_1^{R\mu - \lambda} = 1
\end{equation}
\begin{equation}\label{cover7eqn3}
(\sigma_4\sigma_1)^{\lambda} = 1
\end{equation}
\begin{equation}\label{cover7eqn4}
\sigma_4^{R\mu}(\sigma_5\sigma_4^{-1})^{\lambda} = 1
\end{equation}
\begin{equation}\label{cover7eqn5}
(\sigma_4\sigma_5)^{R\mu}\sigma_5^{-2\lambda} = 1
\end{equation}
\begin{equation}\label{cover7eqn6}
(\sigma_1\sigma_4)^{\lambda} = 1
\end{equation}
Let us also assume for simplicity that $\sigma_5$ commutes with $\sigma_1$.
Equations (~\ref{cover7eqn1}),(~\ref{cover7eqn4}) and (~\ref{cover7eqn5})
 determine an abelian group $A$ of order
 $|(R\mu -2\lambda)^2 - 2\lambda^2|$;
  we must show that $\sigma_4^2$ is non-trivial in $A$.
 $\sigma_4^2$ and $\sigma_5$ generate a subgroup $H$ of $A$ of index
 at most 2.  If $\sigma_4^2 = id$, then $H$ is cyclic of
 order $gcd(|\lambda|, |R\mu - 2\lambda|)$. 
 Then $|(R\mu - 2\lambda)^2 - 2\lambda^2| = |A| \leq 2|H|$
 $ = 2gcd(|\lambda|, |R\mu - 2\lambda|) = 2gcd(|\lambda|, |R|) \leq 2|R|$.
   So if
\begin{equation} \label{sig4eqn}
|(R\mu - 2\lambda)^2 - 2\lambda^2| > 2|R|
\end{equation}
then $\sigma_4^2 \neq id$.
 Therefore, under our hypotheses,
 the relations generate a group which is isomorphic
 to the direct sum of a cyclic group with a hyperbolic triangle group.
  As in the previous cases, we may then find a non-peripheral fixed
 class (in rows 3 and 7), and we are done.

  However, note that if $R = 1$, then Equation (~\ref{sig4eqn})
 is false for all $(\mu, \lambda)$
 satisfying
\begin{eqnarray*}
(\mu+2\lambda)^2 - 2\lambda^2 = 1
\end{eqnarray*}
This is an example of Pell's equation, which has infinitely many solutions,
  and hence (~\ref{sig4eqn}) may be false infinitely often.

\vskip1pc
\textit{Case 5b:} $1/|R\mu - 2\lambda| + 1/|\lambda| < 2/3$
 and $|(R\mu - \lambda)^2 - 2\lambda^2| > 2|R|$.
\vskip1pc

 Let $\sigma_3 = id$, $\sigma_4 = \sigma_2^{-1}$,
 $\sigma_5 = \sigma_1^{-1}$, and $\sigma_7 = \sigma_6^{-1}$ (see Fig.
 9b).
  The conditions for a cover are:
\begin{equation}\label{cover7beqn1}
[\sigma_1,\sigma_2] = 1
\end{equation}
\begin{equation}\label{cover7beqn2}
\sigma_1^{R\mu}(\sigma_2\sigma_1^{-1})^{\lambda} = 1
\end{equation}
\begin{equation}\label{cover7beqn3}
(\sigma_1\sigma_2)^{R\mu}\sigma_2^{-\lambda} = 1
\end{equation}
\begin{equation}\label{cover7beqn4}
(\sigma_6\sigma_1)^{\lambda} = 1
\end{equation}
\begin{equation}\label{cover7beqn5}
\sigma_6^{R\mu - 2\lambda} = 1
\end{equation}
\begin{equation}\label{cover7beqn6}
(\sigma_1\sigma_6)^{\lambda} = 1
\end{equation}
  For simplicity, suppose also that $\sigma_2$ commutes with $\sigma_6$.
Then (~\ref{cover7beqn1}),(~\ref{cover7beqn2}),(~\ref{cover7beqn3})
 determine an abelian group $B$ of order
 $|(R\mu-\lambda)^2 - R\mu\lambda|$.
   If $\sigma_1^2 = id$, then
 $|B| \leq 2gcd(|\lambda|, |R\mu - \lambda|) = 2gcd(|\lambda|, |R|) \leq 2|R|$.
  Therefore, in this case, the group determined by conditions
 (~\ref{cover7beqn1})-(~\ref{cover7beqn6})
 is again the direct product of an abelian group with a hyperbolic
 triangle group, and we may find a non-peripheral fixed class
 in rows 5 and 7.

Note that Case 5a or 5b applies to all but finitely many surgeries
 where $|\lambda| > 1$.  

This concludes the proof of Theorem ~\ref{main}
\end{proof}

\section{Examples}

We begin with the proof of Theorems ~\ref{sister} and ~\ref{fig8}.\\

\begin{lem} ~\label{8lem}
Let $g = D_y^5D_x^{-1}$, $h = D_xD_y$ and $(-1) =  (D_xD_y^{-1}D_x)^2$.
 Then $M_{(-1)h}(\mu,\lambda) \cong M_g(\mu, \lambda - \mu)$, and 
 $M_{h^2}(\mu,\lambda) \cong M_{(-(1)h)^2}(\mu,\lambda+\mu)
 \cong M_{g^2}(\mu,\lambda-\mu)$.
\end{lem}

\begin{proof}
 
Recall that the mapping class group of the once-punctured torus
 is isomorphic to $SL_2(\mathbb{Z})$, under the identifications
 $D_x \rightarrow R=[{1\atop 0} {1\atop 1}]$ and
 $D_y \rightarrow L = [{1\atop 1} {0\atop 1}]$.
  Under these identifications,
 we compute that $h$ has monodromy matrix $[{2 \atop 1} {1 \atop 1}]$,
 $(-1)$ has monodromy matrix $[{-1 \atop 0} {0 \atop -1}]$,
 and $g$ has monodromy matrix $[{1 \atop 5} {-1 \atop -4}]$.
  $h^2$ and $(-h)^2$ have the same monodromy matrix, and hence
 are isotopic.  Therefore $M_{h^2} \cong M_{(-h)^2}$.
 Also, $[{-1\atop -2}{1\atop 1}](-RL)[{-1\atop -2}{1\atop 1}]^{-1}=L^5R^{-1}$,
 so $g$ and $-h$ have conjugate monodromy matrices.
  It follows that $M_{-h} \cong M_g$, and $M_{(-h)^2} \cong M_{g^2}$. 

It remains to determine the effect of these homeomorphisms on the framings.
 Computing the maps on $\pi_1(F)$ gives:
$(-h)^2_{\sharp} = (x^{-1}yxy^{-1})(h^2_{\sharp})(x^{-1}yxy^{-1})^{-1}$.
  Therefore the isotopy which takes $h^2$ to $(-h)^2$ twists $\partial F$
 once in a counter-clockwise manner, so
 the induced bundle homeomorphism sends
 $M_{h^2}(\mu, \lambda)$ to $M_{(-h)^2}(\mu, \lambda+ \mu)$.

 Let $f = D_y^2D_x^{-1}$.
  The bundle homeomorphism induced by conjugation preserves the framing,
 so $M_{-h}(\mu, \lambda)\cong$
 $M_{f(-h)f^{-1}}(\mu, \lambda)$. 
 $f(-h)f^{-1}$ and $g$ have identical
 monodromy matrices, and hence are isotopic.
$g_{\sharp} =$
$(yx^{-1}y^{-1}x)f(-h)f^{-1}_{\sharp}(yx^{-1}y^{-1}x)^{-1}$ 
 so the isotopy from
  $f(-h)f^{-1}$ to $g$ twists
 $\partial F$ once in a clockwise manner.
 The induced bundle homeomorphism sends
 $M_{f(-h)f^{-1}}(\mu, \lambda)$
 to $M_g(\mu, \lambda-\mu)$.
  So $M_{(-1)h}(\mu, \lambda) \cong M_g(\mu, \lambda-\mu)$.
 
  Likewise, $M_{f(-h)^2f^{-1}}(\mu, \lambda) \cong M_{g^2}(\mu, \lambda-2\mu)$.
  Thus $M_{h^2}(\mu, \lambda) \cong M_{(-h)^2}(\mu, \lambda+\mu) \cong$
  $M_{f(-h)^2f^{-1}}(\mu, \lambda + \mu)$
 $\cong M_{g^2}(\mu, \lambda - \mu)$.
\end{proof}

\begin{proof} (of Theorem ~\ref{sister})\\
This is an immediate consequence of Lemma ~\ref{8lem} and
 Theorem ~\ref{main}.
\end{proof}

\begin{proof} (of Theorem ~\ref{fig8})\\
  $M(2\mu,\lambda) \cong M_h(2\mu,\lambda)$, which is double covered
 by $M_{h^2}(\mu, \lambda) \cong M_{g^2}(\mu, \lambda - \mu)$.
  So it is enough to show that $M_{g^2}(\mu, \lambda - \mu)$
 is virtually $\mathbb{Z}$-representable.
  By Theorem ~\ref{main}, we are done unless

\begin{equation*}
1/|-2\mu -(\lambda - \mu)| + 1/|-2\mu - 2(\lambda - \mu)|+1/|\lambda -\mu| \geq 1
\end{equation*}
or, simplifying:
\begin{equation} \label{fig8eqn1}
1/|\mu + \lambda| + 1/|2\lambda| + 1/|\mu - \lambda| \geq 1
\end{equation}

  By [B3], $M(2\mu, \lambda)$ is virtually $\mathbb{Z}$-representable
 if $2\mu$ is divisible by 4; hence we may assume $\mu$ is odd.
 Also, since $gcd(2\mu, \lambda) = 1$, we may assume $\lambda$ is odd,
 and, assuming $(\mu, \lambda) \neq (\pm1,1), \,\, |\lambda| \neq |\mu|$.
  It follows that 
\begin{equation} \label{fig8eqn2}
|\mu - \lambda| \geq 2
\end{equation}
\begin{equation} \label{fig8eqn3}
|\mu + \lambda| \geq 2
\end{equation}
The only simultaneous solutions to inequalities
 (~\ref{fig8eqn1}), (~\ref{fig8eqn2}) and (~\ref{fig8eqn3})
 with $\mu$ and $\lambda$ odd are:
 $(\mu, \lambda) = \pm(-3,1)$ and $\pm(3,1)$.
  So the only possible exceptions to Theorem ~\ref{fig8}
 are $M(-6,1) \cong M(6,1)$ and $M(2,1) \cong M(-2,1)$.
  The virtual $\mathbb{Z}$-representability of these manifolds
 may be verified with either the computer program GAP or Snappea. 
\end{proof}

We now turn to the proof of Theorem ~\ref{fin}.

\begin{proof}(of Theorem ~\ref{fin})\\
 Let $g$ and $h$ be as in the statement of Lemma ~\ref{8lem},
 let $f = h^{18}$, and let $i = D_x^2D_y^{-4}D_xD_y^{-4}D_x$.
  $h^3$ and $i$ both have monodromy
 matrix $[{13\atop{8}}{8 \atop{5}}]$; hence $h^3$ and $i$ are isotopic.
 By arguments similar to those used in the proof of
 Lemma ~\ref{8lem}, we compute that
 $M_{h^3}(\mu, \lambda) \cong M_i(\mu, \mu+\lambda)$.  Therefore
 $M_f(\mu, \lambda) \cong M_{i^6}(\mu,\lambda + 6\mu)$.
  Hence by Theorem ~\ref{main} iii,
 $M_f(\mu, \lambda)$ is virtually $\mathbb{Z}$-representable if
\begin{equation} \label{fineqn1}
 1/|6\mu-\lambda| + 1/|6\mu+\lambda| < 1
\end{equation}

 By Lemma ~\ref{8lem} we have
 $M_f(\mu, \lambda) \cong M_{g^{18}}(\mu, \lambda - 9\mu)$.
 Hence by Theorem ~\ref{main} ii,
 $M_f(\mu, \lambda)$ is virtually $\mathbb{Z}$-representable if
\begin{equation} \label{fineqn2}
1/|9\mu + \lambda| + 1/|2\lambda| + 1/|9\mu - \lambda| < 1
\end{equation}
 
  The only simultaneous solutions to the 
 inequalities ~\ref{fineqn1} and ~\ref{fineqn2}
 have $\mu = 0$.
The proof is completed by noting that $M(0,1)$ has positive first
 Betti number, as it is a torus bundle over $S^1$.

\end{proof}

We remark that the same methods may be applied to many other examples
 of once-punctured torus bundles, to show that all but finitely many
 surgeries are virtually $\mathbb{Z}$-representable.  The idea
 is to start with a monodromy $f$ to which Theorem ~\ref{main}
 i or ii applies.
  Since $L^4$ and $R$ generate a finite-index subgroup of
 $SL_2(\mathbb{Z})$, there exists an integer $\ell$ such
 that $f^{\ell}$ is isotopic to a $g$
 satisfying the hypotheses of Theorem ~\ref{main} iii.
 Usually Theorem ~\ref{main} will then 
 imply that all but finitely many surgeries
 on $M_{f^{\ell}}$ are virtually $\mathbb{Z}$-representable.

\vskip1pc

\section{References}

[B1] M. Baker, \textit{Covers of Dehn fillings on once-punctured torus
 bundles}, Proc. Amer. Math. Soc. \textbf{105} (1989), 747-754.

\vskip1pc

[B2] M. Baker, \textit{Covers of Dehn fillings on once-punctured torus
 bundles}, Proc. Amer. Math. Soc. \textbf{110} (1990), 1099-1108.

\vskip1pc

[B3] M. Baker, \textit{On coverings of figure eight-knot surgeries},
 Pacific J. Math. \textbf{150} (1991), 215-228.

\vskip1pc

[BZ] S. Boyer and X. Zhang, \textit{Virtual Haken 3-manifolds and Dehn
 filling}, preprint.

\vskip1pc

[CJR] M. Culler, W. Jaco, and H. Rubinstein, \textit{Incompressible surfaces
 in once-punctured torus bundles}, Proc. London Math. Soc. (3) \textbf{45}
 (1982), 385-419.

\vskip1pc

[FH] W. Floyd and A. Hatcher, \textit{Incompressible surfaces in punctured
 torus bundles}, Topology and it Applications \textbf{13} (1982), 263-282.

\vskip1pc

[H] J. Hempel, \textit{Coverings of Dehn fillings of surface bundles},
 Topology and its Applications \textbf{24} (1986), 157-170.

\vskip1pc

[KL] S. Kojima and D. Long, \textit{Virtual Betti numbers of some
 hyperbolic 3-manifolds}, A Fete of Topology, Academic Press, 1988.

\vskip1pc

[M] S. Morita, \textit{Finite coverings of punctured torus bundles and
 the first Betti number}, Sci. Papers College Arts Sci. Univ Tokyo
 \textbf{35} (1986), 109-121.

\vskip1pc

[N] A. Nicas, \textit{An infinite family of hyperbolic non-Haken
 3-manifolds with vanishing Whitehead groups}, Math. Proc. Camb. Phil. Soc.,
 \textbf{99} (1986), 239-246.

\vskip1pc

[V] E.B. Vinberg, \textit{Groups defined by periodic paired relations},
 Sbornik: Mathematics, \textbf{188} (1997), 1269-1278.
\\
\\
\\
Joseph D. Masters\\
Mathematics Department\\
University of Texas at Austin 78712\\
masters@math.utexas.edu

\end{document}